\documentclass[krantz1,ChapterTOCs]{krantz-grg} 
\usepackage{fixltx2e,fix-cm}
\usepackage{amssymb}
\usepackage{amsmath}
\usepackage{amsfonts}
\usepackage{amssymb}
\usepackage{amsthm}
\usepackage{subfigure}
\usepackage{graphicx}
\usepackage{color}
\usepackage{enumerate}
\usepackage{lineno}
\usepackage{graphicx}
\usepackage{makeidx}
\usepackage{multicol}
\usepackage{lipsum}
\usepackage{pinlabel}

\newtheorem{theorem}{Theorem}[section] 

\newtheorem{corollary}[theorem]{Corollary}
\newtheorem{algorithm}[theorem]{Algorithm}

\theoremstyle{definition}
\newtheorem{definition}[theorem]{Definition}

\theoremstyle{remark}

\def\O{\mathcal{O}}
\AtBeginDocument{\renewcommand{\bibname}{Resource List}}
\HeadingsChapterSection

\makeindex  

\begin{document}

\chapter*{}

\chapterauthor{Spencer Backman}{Einstein Institute for Mathematics\\
Edmond J. Safra Campus\\
The Hebrew University of Jerusalem\\
Givat Ram. Jerusalem, 9190401, Israel\\      
\texttt{spencer.backman@mail.huji.ac.il}\\}


\chapter{Tutte Polynomial Activities}


\section{Synopsis}

Activities are certain statistics associated to spanning forests and more general objects, which can be used for defining the Tutte polynomial. This chapter is intended to serve as an introduction to activities for graphs and matroids.    We describe  \index{activity}

\begin{itemize}

\item Tutte's original spanning forest activities 

\item Gordon and Traldi's subgraph activities

\item Gessel and Sagan's depth-first search activities

\item Bernardi's embedding activities

\item Gordon and McMahon's generalized subgraph activities

\item Las Vergnas' orientation activities 

\item Etienne and Las Vergnas' activity bipartition

\item Crapo's activity interval decomposition

\item Las Vergnas' active orders

\item Shellability and the algebraic combinatorics of activities

\end{itemize}

\section{Introduction}   
\index{activity} \index{Tutte polynomial} \index{spanning forest} \index{matroid} \index{cut} \index{cycle}
Unlike Whitney's definition of the corank-nullity generating function \index{corank-nullity generating function} $T(G;x+1,y+1)$, Tutte's definition of his now eponymous polynomial $T(G;x,y)$ \index{Tutte polynomial} requires a total order on the edges of which the polynomial is {\it a posteriori} independent.  Tutte presented his definition in terms of \emph{internal and external activities} of maximal spanning forests.  

Although Tutte's original definition may appear somewhat ad hoc upon first inspection, subsequent work by various researchers has demonstrated that activity is a deep combinatorial concept.  In this chapter, we provide an introduction to activities for graphs and matroids. Our primary goal is to survey several notions of activity for graphs which admit expansions of the Tutte polynomial.  Additionally, we describe some fundamental structural theorems, and outline connections to the topological notion of shellability as well as several topics in algebraic combinatorics.  

We use the language of graphs except in sections \ref{active} and \ref{shell} where matroid terminology is employed, although sections \ref{spanningforest}, \ref{activebipartition}, \ref{activesubgraph}, and \ref{unified} apply equally well to matroids, and section \ref{orientation} applies to oriented matroids.

\section{Activities for maximal spanning forests}\label{spanningforest}

We recall Whitney's original definition of the the Tutte polynomial \cite{MR1562461}.  The rank of a subset of edges $A$, written $r(A)$, is the maximum cardinality of a forest contained in $A$. \index{forest}

\begin{definition}\label{EMM:def:ss}
 If $G=(V,E)$ is a graph, then the Tutte polynomial of $G$ is
 \begin{equation}\label{eq:rank_expansion}
     T(G; x, y) = \sum_{A\subseteq E} (x-1)^{r(E)-r(A)}  (y-1)^{|A|-r(A)}.
 \end{equation}
 \end{definition}

 Let $G=(V,E)$ be a graph  and $F$ be  a maximal spanning forest of $G$.
The maximality of $F$ means that if $G$ is connected then $F$ is a spanning tree, and if $G$ is not connected, restricting $F$ to any component of $G$ gives a spanning tree of that component. 

\index{activity} \index{Tutte polynomial} \index{spanning forest} \index{cut!fundamental} \index{cycle!fundamental}
\begin{definition}
Let $F$ be a maximal spanning forest of $G$, $f \in  F$, and $e \in E \setminus F$.   The \emph{fundamental cut associated with $F$ and $e$} \index{cut}\index{cycle} is
 \[  U_F(e) =    \text{\{the edges of the unique cut in $(E \setminus F) \cup e$\}}.  \]

Similarly, the \emph{fundamental cycle associated with $F$ and $e$} is 
\[  Z_F(e) =\text{\{the edges of the unique cycle in $F\cup e$\}}.  \]

\end{definition}

We describe Tutte's activities using fundamental cuts and cycles.

\begin{definition} Let $G = (V,E)$ be a graph with a total order on $E$, and $F$ be a maximal spanning forest of $G$.  We say that an edge $e \in E$ is \index{spanning forest}
\begin{itemize}
\item  \emph{internally active} ($e \in \mathrm I(F)$) with respect to $F$ if $e\in F$ and it is the smallest edge in $U_F(e)$,
\item \emph{externally active} ($e \in \mathrm E(F)$) with respect to $F$ if $e\not\in F$ and it is the smallest edge in $Z_F(e)$.

\end{itemize}
We note that all bridges are internally active and all loops are externally active.
\end{definition}

\begin{figure}[h]
\begin{subfigure}

\centering
\includegraphics[scale=1]{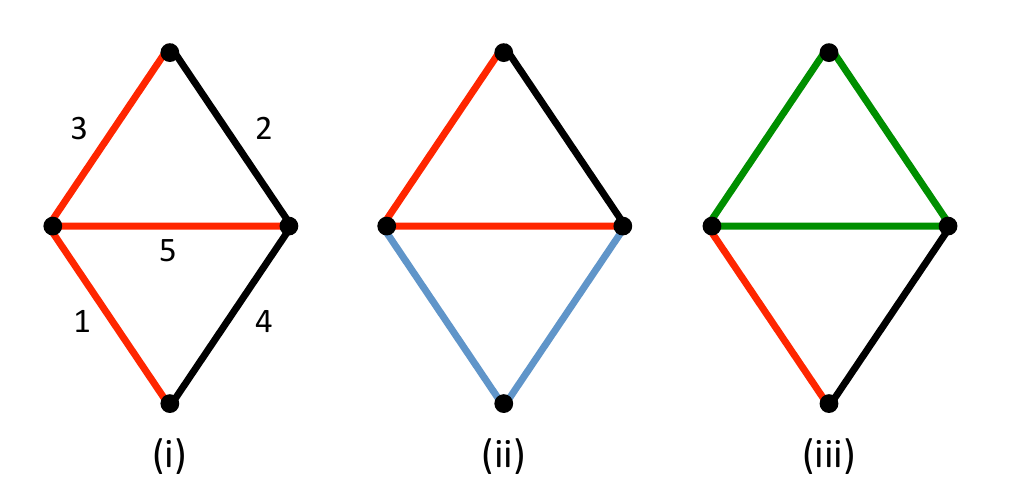}
\caption{(i) A graph with a total order on the edges and a spanning tree $T$ in red. (ii) Edge 1 is the only internally active edge and its fundamental cut is colored blue. (iii) Edge 2 is the only externally active edge and its fundamental cycle is colored green. } 
\end{subfigure}
\end{figure}

 \index{activity}

We can now give the \emph{spanning tree (maximal spanning forest) activities expansion} of the Tutte polynomial.   \index{Tutte polynomial}
\begin{definition}[Tutte \cite{MR0061366}]\label{EMM:d.intext}
If $G=(V,E)$ is a graph with a fixed total order of $E$, then
\begin{equation}\label{EMM:e2}
  T(G;x,y)=\sum\limits_{\substack{F }}x^{|\mathrm{I}(F)|} y^{|\mathrm{E}(F)|},
  \end{equation}
where the sum is over all maximal spanning forests of $G$.   
\end{definition}

Tutte demonstrated that this polynomial is well defined, i.e., it is independent of the total order on the edges. 

\section{Activity bipartition}\label{activebipartition}
\index{activity}\index{activity!bipartition} \index{Tutte polynomial} \index{spanning forest} \index{matroid}\index{flat}\index{convolution}
Etienne and Las Vergnas showed that the activities of a maximal spanning forest induce a canonical bipartition of the edge set of a graph or more generally a matroid. First we recall the definition of a flat.

\begin{definition} Let $\mathcal{F}\subset E$. If there exists no $e \in E \setminus \mathcal{F}$ such that $e$ is contained in a non loop cycle in $e \cup \mathcal{F}$, we say that $\mathcal{F}$ is a \emph{flat} of $G$. A flat is \emph{cyclic} if it is a union of cycles. \index{forest}
\end{definition}

\begin{theorem}(Etienne-Las Vergnas \cite{MR1489076}). Given a maximal spanning forest $F$ of $G$, there exists a unique cyclic flat $\mathcal{F} $ of G such that $\mathcal{F} \cap F$ is a maximal spanning forest of $G|_{\mathcal{F}}$ with no internal activity, and $\mathcal{F}^c \cap F$ is a maximal spanning forest of $G/\mathcal{F} $ with no external activity.
\end{theorem}

As an application of this decomposition, one may obtain a convolution formula for the Tutte polynomial which was independently discovered by Kook, Reiner, and Stanton via incidence algebra methods, and has since been substantially refined and generalized \cite{MR1676189, :ab, MR2718681, MR3678586, MR3887405}.  \index{Tutte polynomial}

\begin{theorem}[Etienne-Las Vergnas \cite{MR1489076}, Kook-Reiner-Stanton \cite{MR1699230}] Given a graph G, then
$$\sum_{\mathcal{F}}T(G;x,y) = T(G/\mathcal{F};x,0)T(G|_{\mathcal{F}};0,y),$$
where the sum is over all cyclic flats of G.
\end{theorem}

\section{Activities for subgraphs}\label{activesubgraph}

\index{subgraph}\index{activity!subgraph}
Gordon and Traldi introduced a notion of activities for arbitrary subgraphs, and used this to provide a 4-variable expansion of the Tutte polynomial which naturally specializes to both Whitney's and Tutte's original expansions.
 \index{activity}

\begin{definition} Let $G = (V,E)$ be a graph with a total order on $E$, and $S \subset E$, then we say that an edge $e$ is
\begin{itemize}
\item  \emph{internally active present} with respect to $S$ $(e \in I(S)\cap S)$ if $e\in S$ and $e$ is the smallest edge in some cut in $S^c \cup e$,
\item  \emph{internally active absent}  with respect to $S$ $(e \in I(S)\cap S^c)$  if $e\notin S$ and $e$ is the smallest edge in some cut in $S^c$,
\item \emph{externally active present}  with respect to $S$ $(e \in  L(S)\cap S)$ if $e\in S$ and $e$ is the smallest edge in some cycle in $S$,
\item \emph{externally active absent}  with respect to $S$ $(e \in \mathrm L(S)\cap S^c)$ if $e\notin S$ and $e$ is the smallest edge in some cycle in $S \cup e$.

\end{itemize}

\end{definition}

\begin{theorem}[Gordon-Traldi \cite{MR1080623}]\label{EMM:d.intext} \index{Tutte polynomial}
If $G$ is a graph with a fixed total order of $E$, then
\begin{equation}\label{EMM:e2}
  T(G;x+w,y +z)=\sum\limits_{\substack{S \subset E}} x^{|{I(S)\cap S}|}w^{| I(S)\cap S^c|} y^{| L(S)\cap S|}z^{|L(S)\cap S^c|}.
  \end{equation}

\end{theorem}

By setting $x= 1$ and $z=1$, we recover Whitney's definition, and by setting $w= 0$ and $y=0$, we recover Tutte's definition.  

 While Gordon and Traldi's expansion is proven recursively via deletion-contraction, from which they obtain more general formulae, the 4-variable expansion is equivalent to an earlier theorem of Crapo.   \index{forest}

\begin{theorem}[Crapo \cite{MR0262095}]

Let $P(E)$ be the Boolean lattice of subgraphs of $G$ ordered by containment.  Given a spanning forest $F$, define an interval in this lattice $[ F \setminus I(F), F \cup E(F)]$.  Then 
$$P(E) = \sqcup_{F} [ F \setminus I(F), F \cup E(F)]$$
where the disjoint union is over all maximal spanning forests.
\end{theorem}

\begin{figure}[h]
\begin{subfigure}

\centering
\includegraphics[scale=.9]{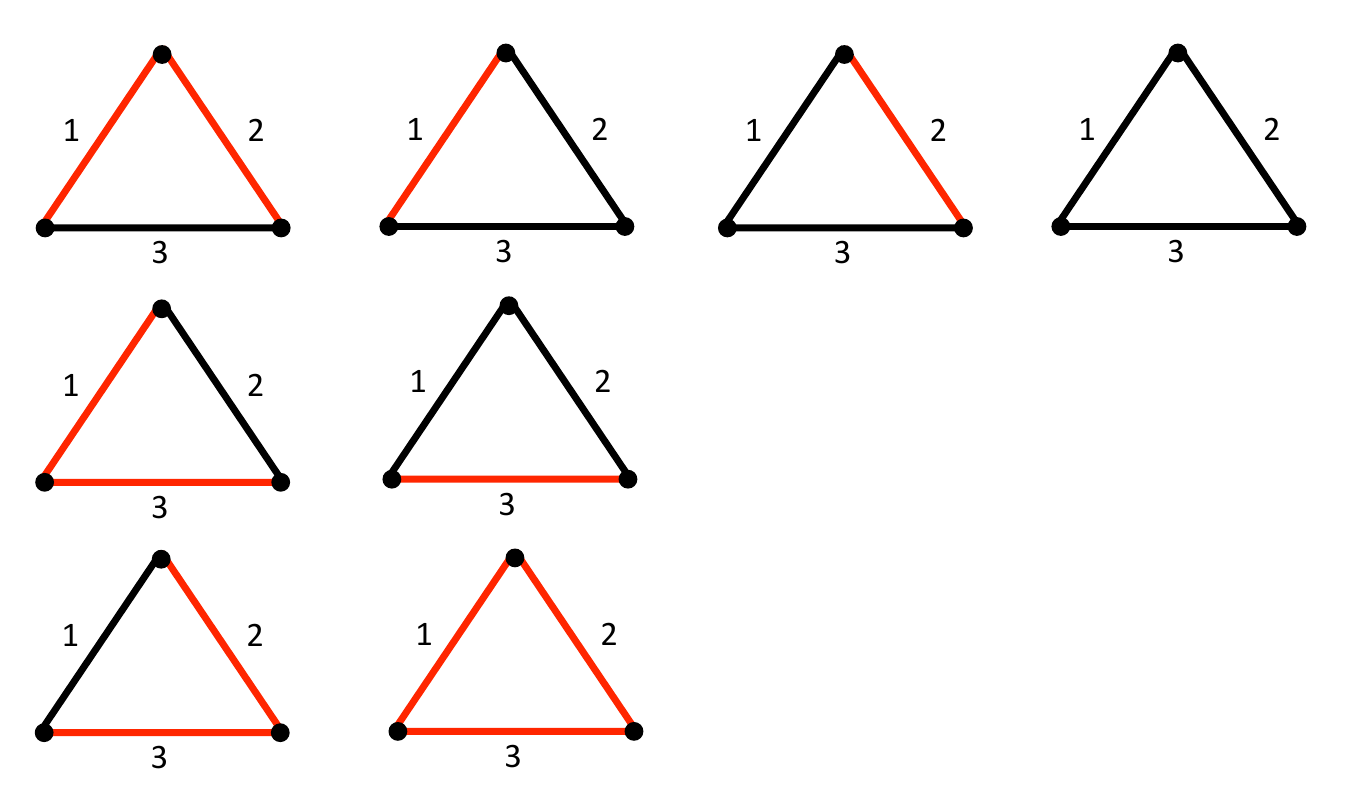}
\caption{A Crapo decomposition of the subgraphs of $K_3$.  Each row collects the subgraphs corresponding to the spanning tree on the left.}
\end{subfigure}
\end{figure}

  In the following sections \ref{DFS}, \ref{combmaps}, and \ref{unified} we will describe other notions of activity for maximal spanning forests which employ input data different from a total order on the edges, but still allow for expansions of the Tutte polynomial and an analogue of Crapo's interval decomposition.
  \index{activity}\index{depth-first search}\index{Crapo interval}\index{spanning tree}\index{activity!depth-first search}
  \section{Depth-first search external activity}\label{DFS}
  
  Gessel and Sagan introduced a notion of external activity for maximal spanning forests based on depth-first search.  For simplicity sake, we will assume that our graph is connected and has no parallel edges.  In what follows, we view a tree rooted at a vertex $q$ to be oriented so that every vertex is reachable from $q$ by a directed path.    \index{activity}
  
  \begin{definition}[Gessel-Sagan \cite{MR1392494}]
  Let $<$ be a total order on the vertices of $G$, and $F$ be a spanning tree of $G$ rooted at the smallest vertex $q$.  Let $e = (u,v)$ be an edge of $G \setminus F$.  We say that $e$ is {\emph{depth-first search externally active (DFS externally active)}}, and write $e\in E_{DFS}(F)$, if either $u=v$, or $(u,w)$ is an oriented edge in $F$ belonging to the unique cycle in $F \cup (u,v)$, and $w>v$.
  
  \end{definition}
  
  \begin{figure}[h]
  \begin{subfigure}

\centering
\includegraphics[scale=.8]{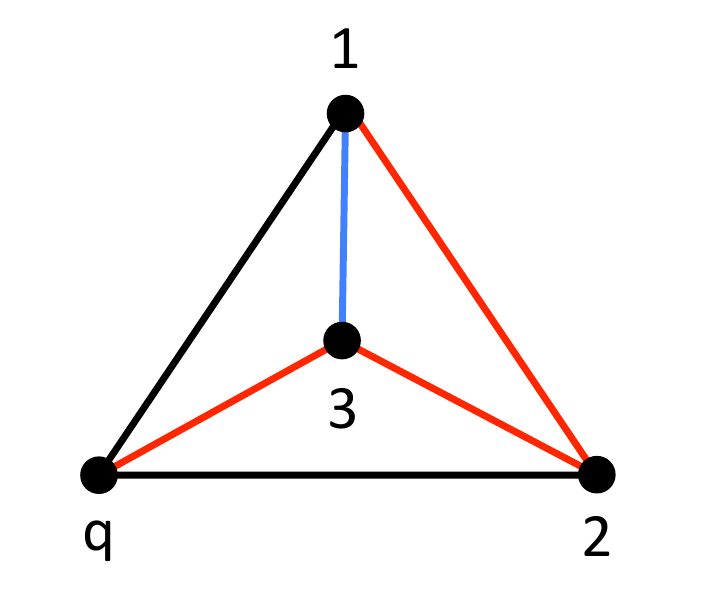}
\caption{$K_4$ with a root $q$ and a total order on its vertices.  A spanning tree in red and a DFS externally active edge in blue.}
\end{subfigure}
\end{figure}

    The name DFS externally active is justified by the following observation: given a spanning subgraph of $G$, we can produce a spanning forest $F$ by performing a DFS search which favors larger labeled vertices.  Then $(u,v)$ is DFS externally active if when we apply DFS search to the graph $F \cup (u,v)$, we obtain $F$.  Gessel and Sagan showed that DFS external activity when combined with Tutte's notion of internal activity allows for an expansion of the Tutte polynomial.\index{spanning forest}\index{Tutte polynomial}
    
    \begin{theorem}[Gessel-Sagan \cite{MR1392494}]\index{Tutte polynomial} \index{forest}
    If $G$ is a connected graph with a total order of its vertices, then
    \begin{equation}\label{GS}
  T(G;x,y)=\sum\limits_{\substack{T}} x^{|{I(T)}|}y^{|{E_{DFS}(T)}|},
  \end{equation}
  where the sum is over all spanning trees $T$ of $G$.
  \end{theorem}

     In the same article, Gessel and Sagan produced a notion of Nearest Neighbor First activity, which will not be reviewed here.

\index{activity!depth-first search}\index{depth-first search}\index{spanning tree}\index{map}\index{embedding}\index{half-edge}\index{activity!embedding}\index{spanning tree}

\section{Activities via combinatorial maps}\label{combmaps}
 \index{activity}
Bernardi proposed a notion of activity induced by a rooted combinatorial map, which is essentially an embedding of a graph $G=(V,E)$ with a distinguished half-edge $h$ into an orientable surface.  In what follows, we assume that $G$ is connected and loopless, although these restrictions are not essential.

Informally, given a spanning tree $T$ of $G$, we can use a rooted combinatorial map to tour the edges of $G$ by starting at $h$ and traveling counterclockwise around the outside of $T$.  We then declare an edge $e \notin T$ to be externally active if it is the first edge in its associated fundamental cycle which we meet during the tour.  We similarly say an edge $e \in T$ is internally active if it is the first edge in its associated fundamental cut which we meet during the tour.

We now describe Cori's maps \cite{MR0404045} and Bernardi's activities more formally.  
We define a \emph{half-edge} to be an edge and an incident vertex, e.g. if $e \in E$ with  $e = (u,v)$  such that $u,v \in V$, then $(e,u)$ and $(e,v)$ are the two associated half-edges.  Let $\sigma$ be a permutation of the half-edges of $G$ such that $\sigma((e,u))=(e',u)$ some other half-edge incident to $u$, and for any two half-edges $(e,u)$ and $(e',u)$ with the same endpoint, there exists some $k >0$ such that  $\sigma^k((e,u))=(e',u)$.  Let $h$ be a distinguished half-edge.  We define a \emph{rooted combinatorial map} to be a triple $(G,\sigma,h)$.  Let $\alpha$ be the involution on the set of half-edges such that for all $e = (u,v) \in E$, $\alpha((e,u)) = (e,v)$.  Given a half-edge $i$ and a spanning tree $T$, we define the \emph{motion operator} \[
    t(i)= 
\begin{cases}
   \sigma(i) &\text{if } i \notin T\\ 
  \sigma \circ \alpha(i)  & \text{if } i \in T           
\end{cases}
\]

It is easy to check that the iterated motion operator defines a tour of the half-edges of $G$.  This tour induces a total order $<_{T}$ on the edges given by the first time one of its half-edges is visited.  We use this total order to define internally and externally active edges. \index{forest}

 \begin{figure}[!ht]
\begin{subfigure}

\centering
 \includegraphics[scale=.35]{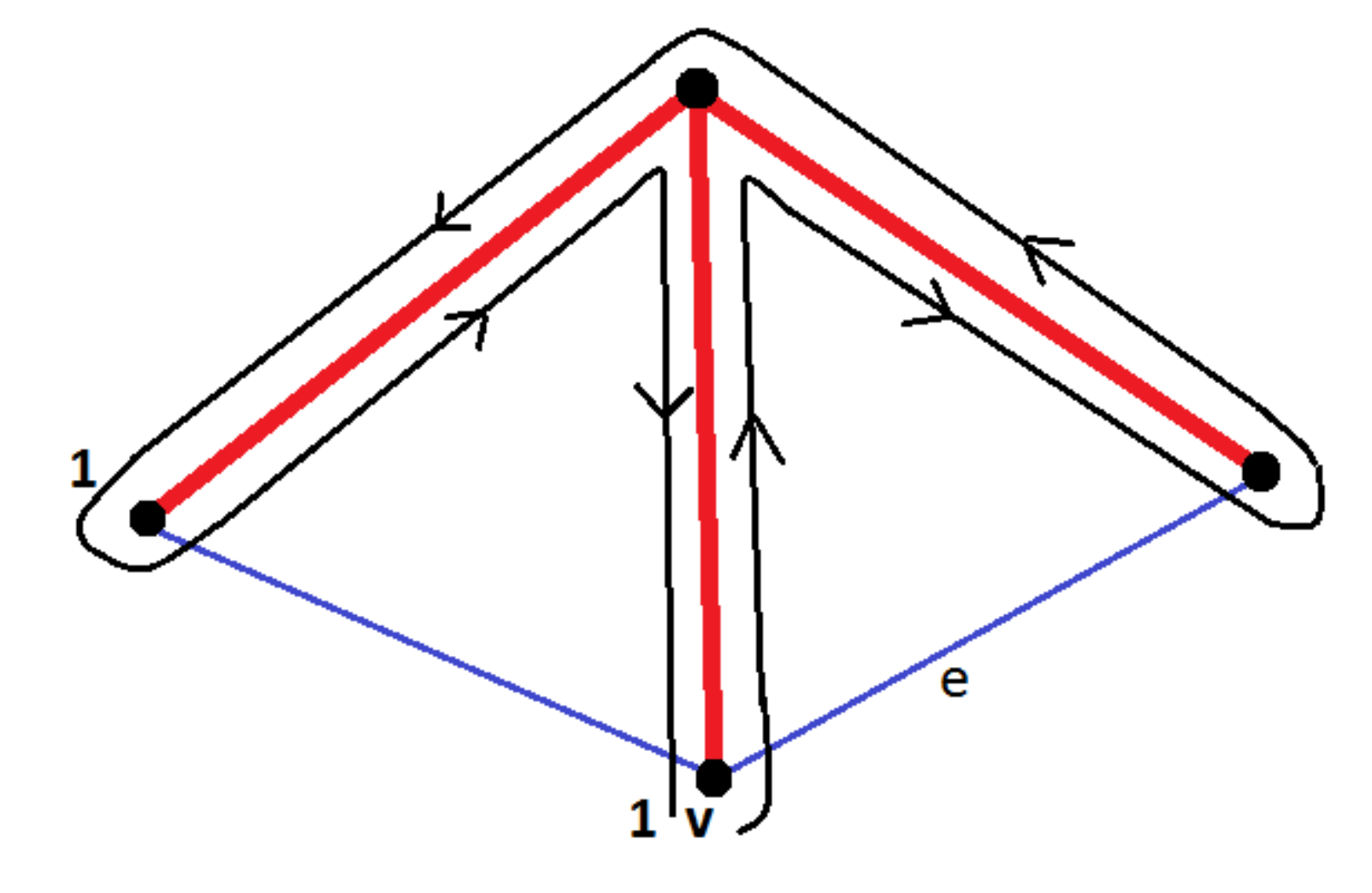}
\end{subfigure}
\begin{subfigure}
\centering
 \includegraphics[scale=.35]{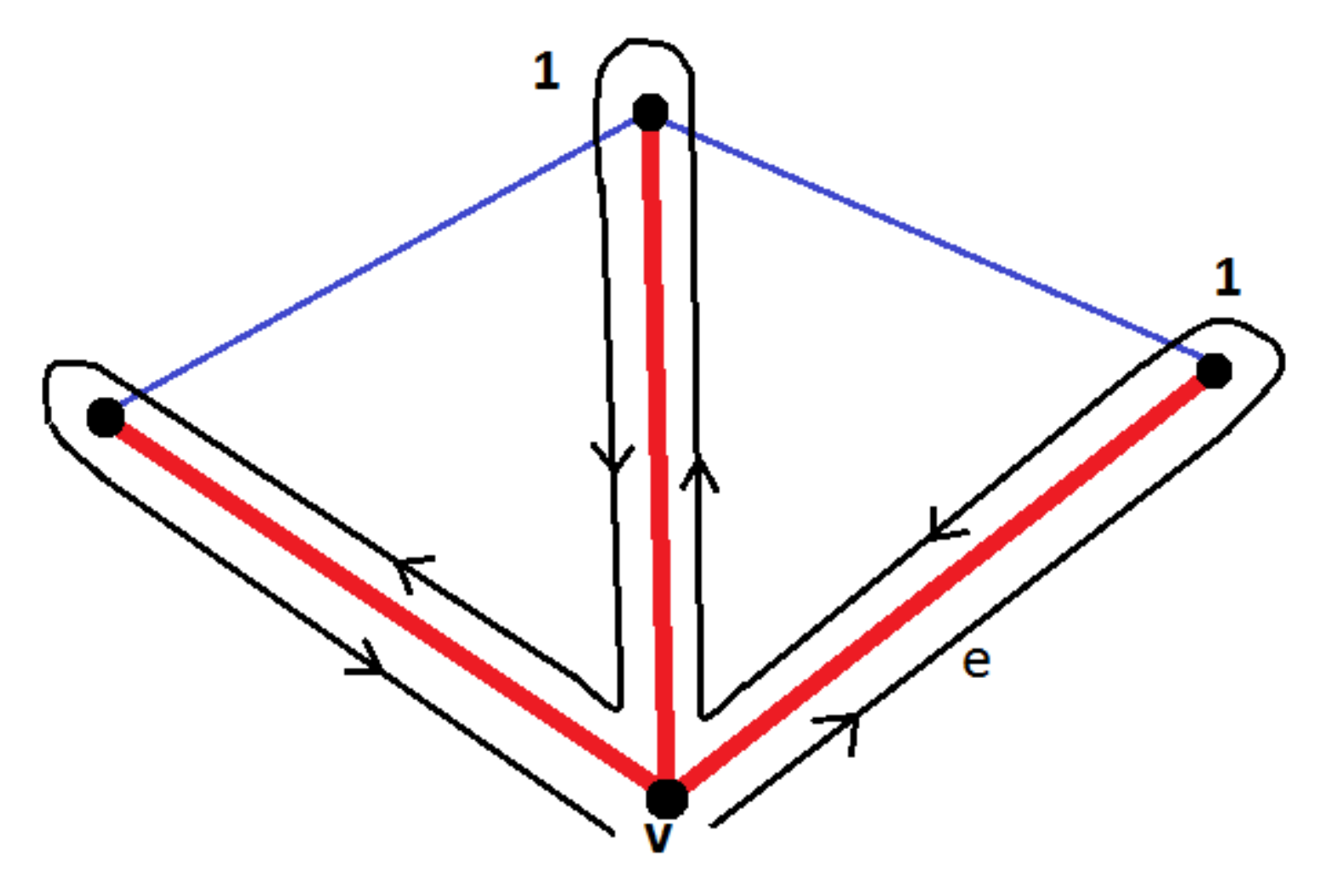}
\caption{Tours of two different spanning trees induced by the motion operator associated to the same rooted combinatorial map.}
\label{bernardi}
\end{subfigure}

\end{figure}
 \index{activity}

\begin{definition}[Bernardi \cite{MR2428901}]  Let $(G,\sigma, h)$ be a rooted combinatorial map, $T$ a spanning tree of $G$, and $e$ an edge of $G$.

\begin{itemize}

\item The edge $e$ is \emph{embedding-internally active} $(e \in I_B(T))$ if $e \in T$ and $e$ is the minimum edge with respect to $<_{T}$ in its associated fundamental cut.

\item The edge $e$ is \emph{embedding-externally active} $(e \in E_B(T))$ if $e \notin T$ and $e$ is the minimum edge with respect to $<_{T}$ in its associated fundamental cycle.

\end{itemize}

\end{definition}

Bernardi showed that this definition of activity admits an expansion of the Tutte polynomial.

    \begin{theorem}[Bernardi \cite{MR2428901}]\index{Tutte polynomial}
    Let $(G,\sigma,h)$ be a rooted combinatorial map, then
    
        \begin{equation}\label{GS}
  T(G;x,y)=\sum\limits_{\substack{T}} x^{|{I_B(T)}|}y^{|{E_B(T)}|},
  \end{equation}
  where the sum is over all spanning trees.
\end{theorem}

We remark that Courtiel recently introduced a different notion of activity via combinatorial maps which he calls the ``blossom activity" \cite{courtiel:hal-01088871}.
\index{Tutte polynomial}

\section{Unified activities for subgraphs via decision trees}\label{unified}
The Gordon-Traldi activities were further generalized by Gordon-McMahon \cite{MR1425948} in a way which also applies to greedoids.  The Gordon-McMahon notion of activity was rediscovered by Courtiel \cite{courtiel:hal-01088871} who proved that it allows for a unification of all of the aforementioned notions of activity.  We describe these activities following Courtiel and using the language of decision trees. 

\index{activity}\index{activity!blossom}\index{activity!unified}\index{decision tree}\index{spanning tree}\index{binary tree}

\begin{definition} \index{activity}
A \emph{decision tree} $D$ for $G$ consists of a perfect binary tree, i.e., a rooted tree such that all non-leaf nodes have two descendants, and a labeling of each node of the tree by elements of $E$ such that the labels along any particular branch give a permutation of $E$.

\end{definition}

  \begin{figure}[h]

\begin{subfigure}

\centering

\includegraphics[scale=.9]{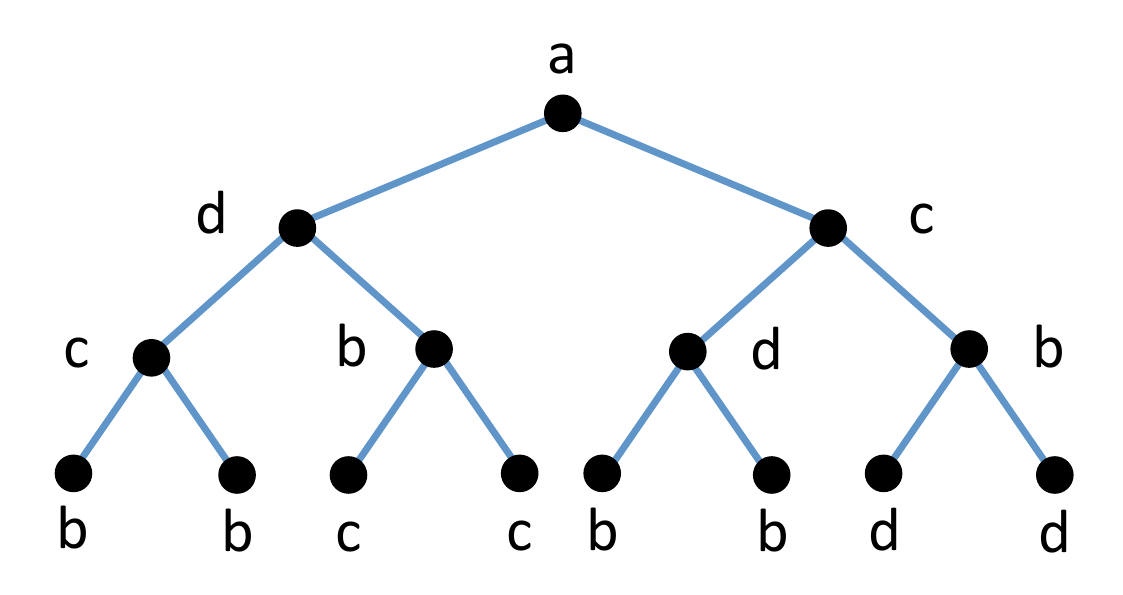}
\caption{An example of a decision tree for a graph on 4 edges.}
\end{subfigure}

\end{figure}

Given a decision tree $D$ and a subgraph $S \subset E$, we can use $D$ to partition $S$ into four sets: $I(S)$, $L(S)$, $S_E$, $S_I$.  We describe the recursive algorithm for producing this partition informally and refer the reader to \cite{courtiel:hal-01088871} for a pseudocode presentation.  

\begin{algorithm}{\text Recursive generalized activities algorithm}

\

Initialize with $X = E, I(S) = L(S) = S_I = S_E = \emptyset$, and $e$ corresponding to the label of the root of $D$.  While $X \neq \emptyset$, do the following:
\begin{itemize}

\item If $e$ is a bridge, add $e$ to $I(S)$, contract $e$ in $X$, and move to the right descendant of $e$ in $D$.
\item If $e \in S $ is neither a bridge nor a loop, add $e$ to $S_I$, contract $e$ in $X$, and move to the right descendant of $e$ in $D$.
\item If $e$ is a loop, add $e$ to $L(S)$, delete $e$ in $X$, and move to the left descendant of $e$ in $D$.
\item  If $e \notin S$ is neither a bridge nor a loop, add $e$ to $S_E$, delete $e$ in $X$, and move to the left descendant of $e$ in $D$.

After we move to a descendant, we update $e$ to be the label of the new node and recurse. 
\end{itemize}

\end{algorithm}
 \index{activity}

\begin{theorem}[Gordon-McMahon \cite{MR1425948}]\index{Tutte polynomial}
Let $G$ be a graph, $D$ a decision tree for $G$, and $I(S)$ and $L(S)$ as above, then the Tutte polynomial has the following expansion

$$T(G;x+w, y+z) = \sum_{S \subset E} x^{|S \cap I(S)|}w^{|S^c \cap I(S)|}y^{|S^c \cap L(S)|}z^{|S \cap L(S)|}.$$

\end{theorem}
\index{activity!orientation}\index{algorithm}\index{algorithm!recursive}\index{Tutte polynomial}\index{orientation}\index{cycle!directed}\index{cut!directed}
\section{Orientation activities}\label{orientation} \index{activity}
A famous result of Stanley states that $T(G;2,0)$ (equivalently, the chromatic polynomial evaluated at $-1$) counts the number of acyclic orientations of $G$ \cite{MR0317988}.  This result was generalized to hyperplane arrangements by Zaslavsky \cite{MR0357135}, and to oriented matroids by Las Vergnas \cite{MR586435}.

Las Vergnas \cite{MR776814} introduced a notion of orientation activities, which parallels those of subsets, and allows for an orientation expansion of the Tutte polynomial which recovers Stanley's result.  He later introduced refined orientation activities \cite{Vergnas:aa}, which we now describe.  Similar to the way that Tutte's activities are defined in terms of fundamental cuts and cycles, Las Vergnas' orientation activities are defined in terms of directed cuts and cycles.

\begin{definition}
Let $\O$ be an orientation of the edges of $G$.  Let $Z$ be a cycle in $G$, and $U = (X,X^c)$ be a cut in $G$.  We say that $Z$ is a \emph{directed cycle} if we can walk around the cycle traveling in the direction of the edge orientations.  We similarly define $U$ to be a \emph{directed cut} if all of its edges are, without loss of generality, oriented from $X$ to $X^c$.

\end{definition}

We now use directed cuts and cycles to introduce notions of orientations activities.
\begin{definition}[Las Vergnas \cite{Vergnas:aa}]\label{oract} \index{activity}\index{orientation}\index{activity!orientation}\index{orientation!reference}
 Let $G = (V,E)$ be a graph with a total order on $E$, and  $\O_{ref}$ a reference orientation of the edges of $G$.  If $\O$ is an orientation of $G$ and $e \in E$, then we say $e$ is
\begin{itemize}
\item  \emph{positive cut active} ($e \in I(\O)^+$)  if $e$ is the smallest edge in some directed cut and is oriented in agreement with $\O_{ref}$,
\item  \emph{negative cut active} ($e \in I(\O)^-$) if $e$ is the smallest edge in some directed cut and is oriented in disagreement with $\O_{ref}$,
\item \emph{positive cycle active} ($e \in E(\O)^+$) if $e$ is the smallest edge in some directed cycle and is oriented in agreement with $\O_{ref}$,
\item \emph{negative cycle active} ($e \in E(\O)^-$) if $e$ is the smallest edge in some directed cycle and is oriented in disagreement with $\O_{ref}$.

\end{itemize}

\end{definition}

  \begin{figure}[h]
  \begin{subfigure}
  
\centering
\includegraphics[scale=1.2]{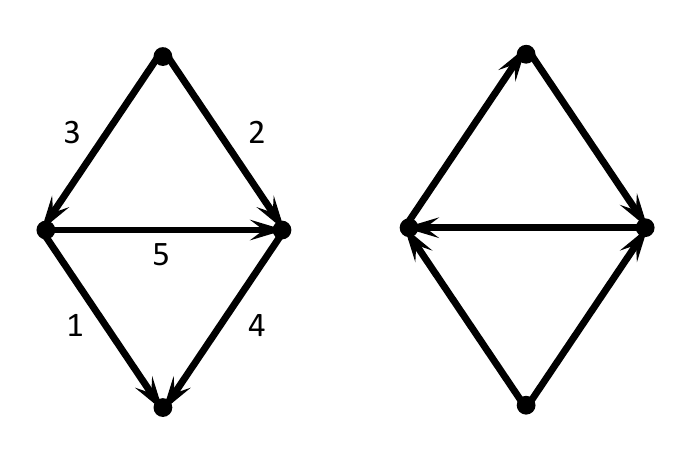}
\caption{Left: A total order and reference orientation of the edges of a graph $G$.  Right: An orientation $\O$ of $G$ with $I(\O)^+ = \emptyset , I(\O)^- = \{ 1\}, L(\O)^+ =\{2\}, {\rm and \, } L(\O)^- = \emptyset$.}

\end{subfigure}
\end{figure}

\begin{theorem}[Las Vergnas \cite{Vergnas:aa}]\label{EMM:d.intext}\index{Tutte polynomial}
Let $G$ be a graph with a fixed total order and reference orientation of $E$, and $I(\O)^+, I(\O)^-, L(\O)^+,$ and $L(\O)^-$ as in Definition \ref{oract}, then
\begin{equation}\label{EMM:e2}
  T(G;x +w,y +z)=\sum\limits_{\substack{\O }} x^{|I(\O)^+|}w^{|I(\O)^-|} y^{|E(\O)^+|}z^{|E(\O)^-|},
  \end{equation}
 where the sum is over all orientations of $E$.

\end{theorem}

Las Vergnas orientation expansion holds for all oriented matroids.  By specializing variables $x=w = u/2$ and $y = z = v/2$, we recover Las Vergnas earlier expansion which does not make use of a reference orientation. 
   
\begin{corollary}[Las Vergnas \cite{MR776814}]
Let $G$ be a graph with a fixed total order $E$, and let $I(\O)$ and $L(\O)$ to be the set of edges which are minimum is some directed cut or cycle, respectively, then
\begin{equation}\label{EMM:e2}
  T(G;x,y)=\sum\limits_{\substack{\O}} \left({u \over 2}\right)^{|I(\O)|}\left({v\over 2}\right)^{|L(\O)|},
  \end{equation}
where the sum is over all orientations of $E$.
\end{corollary}

\index{activity}\index{orientation}\index{activity!orientation}\index{active bijection}\index{fourientation}\index{subgraph}\index{Tutte polynomial}\index{active orders}\index{pivoting}\index{matroid}\index{activity!fourientation}

We remark that Berman \cite{MR0469810} was the first to propose an orientation expansion of the Tutte polynomial, although his definition was not correct.  There are natural notions of orientation activity classes which parallel Crapo's subset intervals.   The refined orientation expansion of the Tutte polynomial also follows as a direct consequence of ``the active bijection" of Gioan and Las Vergnas which gives a bijection between orientation activity classes and Crapo subset intervals, which respects the four different activities \cite{MR2552669, MR3895683}.  See Gioan's chapter for further details.  \index{activity}

A \emph{fourientation} of a graph is a choice for each edge of the graph whether to orient that edge in either direction, leave it unoriented, or biorient it.  One may naturally view fourientations as a mixture of orientations and subgraphs where absent and present edges correspond to unoriented and bioriented edges, respectively.  Backman, Hopkins, and Traldi \cite{MR3707217} introduced notions of activities for fourientations which provide a common refinement of the Gordon and Traldi subgraph activities and Las Vergnas' orientation activities.

\section{Active orders}\label{active}
 \index{activity}
For deepening our understanding of activities, Las Vergnas introduced three \emph{active orders} on the bases of a matroid.  We describe these partial orders here and Las Vergnas' key result that they induce lattice structures on bases.

Let $B_1$ and $B_2$ be bases of a matroid $M$ with fixed order on its ground set.  We say that $B_1$ is obtained from $B_2$ by an \emph{externally active pivoting} if $B_2 = B_1 \setminus e \cup f$, where $e$ is the minimum element in $Z_{B_1}(f)$, and we write $B_1 \leftarrow_M B_2$.  Dually, we say that $B_1$ is obtained from $B_2$ by an \emph{internally active pivoting} if $B_1 = B_2 \setminus e \cup f$, where $e$ is the smallest element in $U_{B_1}(f)$, and we write $B_1 \leftarrow *_M B_2$.  We let $<_{Ext}$ and $<_{Int}$ denote the partial orders on the bases obtained by taking the transitive closures of the relations $ \leftarrow_M$ and $ \leftarrow *_M$, respectively.  We refer to $<_{Ext}$ as the \emph{external order}, and $<_{Int}$ as the \emph{internal order}.

Las Vergnas also defined the following join of the external and internal orders.  Let $B_1$ and $B_2$ be bases, then we say that $B_1 <_{Ext/Int} B_2$ if there exists bases $C_1, \dots , C_k$ such that $B_1 = C_1$, $B_2 = C_k$, and for each $i$ either $C_i  \leftarrow_M C_{i+1}$ or $C_i  \leftarrow *_M C_{i+1}$.

Recall that a lattice is a poset such that every pair of elements have a unique join and meet.

\begin{theorem}[Las Vergnas \cite{MR1845495}]
Let $\mathcal{B} (M)$ be the set of bases of a matroid $M$ with a total order on its ground set, then the posets \mbox{$(\mathcal{B} (M) \cup \{\mathbf{0}\}, <_{Ext})$}, \mbox{$(\mathcal{B} (M) \cup \{\mathbf{1}\}, <_{Int})$}, and $(\mathcal{B} (M), <_{Ext/Int})$ are lattices.

\end{theorem}
\index{active orders}\index{pivoting}\index{matroid}\index{bases}\index{lattice}\index{lattice!atomistic}\index{lattice!supersolvable}\index{lattice!distributive}\index{shellability}\index{simplicial complex}\index{facet}\index{face}\index{complex!independence}\index{complex!no broken circuit}\index{$f$-polynomial}\index{$h$-polynomial}\index{$f$-vector}\index{$h$-vector}\index{Tutte polynomial}
Las Vergnas observed that the lattice associated to the external order is not distributive, although it is atomistic.  This appears to have been remedied in the recent PhD thesis of Gillespie \cite{:aa} where it is shown that by extending the external order to all independent sets, a supersolvable join-distributive lattice is obtained.

\section{Shellability and activity}\label{shell}
 \index{activity}
There are important connections between activity and combinatorial topology.  We refer the reader to Bj\"orner \cite{MR1165544} for an excellent introduction to this topic.  Let $[n]$ denote the finite set $\{1, \dots , n\}$.  An \emph{abstract simplicial complex} (often just called a simplicial complex) $\Delta$ on $n$ elements is a collection of subsets (faces) of $[n]$ which is closed under taking subsets.

Informally, a simplicial complex $\Delta$ is \emph{shellable} if there is an ordering of the maximal faces $(\emph{facets})$ of $\Delta$ so that each facet can be added to the previous ones by glueing along codimension 1 faces.

The \emph{f-polynomial} of a simplicial complex $\Delta$ is $f_{\Delta}(x) = \sum_{i=0}^n f_ix^{d-i}$ where $f_i$ is the number of faces of $\Delta$ of size $i$.  The \emph{h-polynomial} of $\Delta$ is the $h_{\Delta}(x) = f_{\Delta}(x-1)$.  The $f$-vector and $h$-vector of $\Delta$ are the vectors whose entries are the coefficients of the $f$-polynomial and $h$-polynomial, respectively.  A shellable complex is homotopy equivalent to a wedge of spheres, and its $h$-vector is nonnegative.

The following complexes are shellable, and the proofs of shellability are related to activities.

\begin{itemize}
\item The \emph{independence complex} $(IN(M))$ whose faces are the independent sets of a matroid $M$ \cite{MR2627467, MR593648}.  Its $h$-polynomial is $T(M;x,1)$.  Matroids are characterized by the fact that the $IN(M)$ are the pure simplicial complexes which are lexicographically shellable with respect to any order on the ground set \cite{MR1165544}.

\item The \emph{no broken circuit complex} $(NBC(M))$ whose faces are the independent sets of a matroid with no external activity in the sense of Tutte \cite{MR2627467, MR593648}.  Its $h$-polynomial is $T(M;x,0)$.

\item The \emph{external activity complex} defined on $E(M)\times E(M)$, whose facets are given by $B \cup L(B) \times B\cup (B \cup L(B))^c$, where $B$ ranges over the bases of $M$ \cite{MR3558045}.  The shelling makes use of Las Vergnas' order $<_{Ext/Int}$.  It has the same $h$-polynomial as $IN(M)$ \cite{MR3558045}.

\item The \emph{order complex of IN(M)} modulo Las Vergnas' external active order $<_{Ext}$ \cite{:aa}.

\item The \emph{order complex of the lattice of flats} \cite{MR570784}.  It has Euler characteristic $T(M;1,0)$.

\end{itemize}
\index{Tutte polynomial}
These complexes admit many interesting connections with algebraic and geometric combinatorics.  We briefly mention a few.  Orlik and Solomon \cite{MR558866} introduced a certain graded algebra which is isomorphic to the cohomology ring of the complement of a complex hyperplane arrangement, and showed that monomials corresponding to faces of $NBC(M)$ give a basis for this algebra; see Falk and Kung's chapter for an introduction to these algebras. \index{activity}\index{shellability}\index{simplicial complex}\index{facet}\index{face}\index{complex!independence}\index{complex!no broken circuit}\index{$f$-polynomial}\index{$h$-polynomial}\index{complex!order}\index{algebra!Orlik-Solomon}\index{tropical geometry}\index{Bergman fan}\index{log-concavity}\index{unimodular}\index{$h$-vector}\index{$f$-vector}\index{$O$-sequence}\index{flat}\index{lattice}\index{monomial}

The external activity complex was introduced by Ardila and Boocher in their investigation of commutative algebraic aspects of the closure of a linear space in a product of projective lines \cite {MR3439307}.  The ideals they consider are homogenizations of ones considered earlier by Proudfoot and Speyer \cite{MR2246531} and Terao\cite{MR1899865}.  A slight variation of these objects play an important role in Huh and Wang's proof of the Dowling-Wilson conjecture for realizable matroids \cite{MR3733101}.

In the field of tropical geometry, which is often referred to a as a piecewise linear version of algebraic geometry, Bergman fans are certain balanced polyhedral fans which provide a new and exciting framework for studying matroids.  Ardila and Klivans \cite{MR2185977} showed that they are unimodularily triangulated by the fan over the order complex of the lattice of flats.  This relationship has lead authors to uncover interesting connections between algebraic geometry and matroids, most notably the proof of the Heron-Rota-Welsh conjecture that the $f$-vector of $NBC(M)$ is log-concave by Adirprasito, Huh, and Katz \cite{MR3862944} building on earlier works of Huh \cite{MR2904577}, and Huh and Katz \cite{MR2983081}.  Recently, Fink, Speyer, and Woo introduced an \emph{extended NBC complex} in order to shed some light on these different manifestations of the $f$-polynomial of $NBC(M)$ \cite{fink}.

Finally, a major open question in combinatorial commutative algebra is Stanley's 1977 conjecture \cite{MR0572989} that the $h$-vector of $IN(M)$ is a pure $O$-sequence.  Merino settled this conjecture in the case of cographic matroids by application of chip-firing \cite{MR1888777}; see Merino's chapter for further details.

\section{Acknowledgements}
Many thanks to the anonymous referee for providing very helpful feedback on an earlier draft of this chapter, and to Sam Hopkins for pointing out several typographical errors.  Additional thanks to Matt Baker for generously sharing Figure \ref{bernardi}.

\bibliographystyle{plain}
\bibliography{tutteactivitiesfinal}

\printindex

\end{document}